\documentclass{article}
\usepackage{amsthm}
\usepackage{fancyhdr}
\usepackage{amsmath}
\usepackage{indentfirst}
\usepackage{titlesec}
\usepackage{amsfonts}
\usepackage{amssymb}

\newcommand{\cA}{\mathcal A}
\newcommand{\cB}{\mathcal B}

\newcommand{\cF}{\mathcal F}

\newcommand{\cS}{\mathcal S}

\newcommand{\N}{\mathbb N}
\newcommand{\R}{\mathbb R}
\newcommand{\Z}{\mathbb Z}
\newcommand{\cov}{\textrm{\rm Cov}}
\newcommand{\corr}{\textrm{\rm Corr}}

\newcommand{\zhfb}{\hfil\break}
\newcommand{\zeps}{\varepsilon}
\newcommand{\bbP}{\mathbb P}
\newcommand{\var}{\textrm{\rm Var}}

 \newtheorem{theorem}{Theorem}[section]

 \newtheorem{lemma}[theorem]{Lemma}

 \theoremstyle{definition}
 \newtheorem{definition}[theorem]{Definition}

\newtheorem{remark}[theorem]{Remark}

 \newtheorem{notations}[theorem]{Notations}

 \newtheorem{setting}[theorem]{Setting}

 \numberwithin{equation}{section}

 \numberwithin{theorem}{section}

\begin{document}

\title{\textbf{On some possible combinations of mixing rates for strictly stationary, reversible Markov chains}}

\newcommand{\orcidauthorA}{0000-0000-000-000X} 

\author{Richard C.\ Bradley \\
Department of Mathematics, 
Indiana University, Bloomington, IN 47405, U.S.A.\\ bradleyr@indiana.edu}


\maketitle

\begin{abstract} 
A class of examples is constructed to show that for strictly stationary Markov chains that are reversible, the simultaneous mixing rates for the $\rho$-mixing and strong mixing ($\alpha$-mixing) conditions can be fairly arbitrary, within certain unavoidable tight restrictions.
The examples constructed here have the added property that the mixing rate for the absolute regularity ($\beta$-mixing) condition is within a constant factor of that for strong mixing.     
\end{abstract}

\hfil\break 
\null \hskip 1.2 in 
strong mixing; absolute regularity; $\rho$-mixing
\medskip





\section{Introduction}
\label{sc1}

     Theory developed in the papers by 
Roberts and Rosenthal \cite{ref-journal-RR} 
and Roberts and Tweedie \cite{ref-journal-RT} together produced, among other things, a ``key result'' giving, for strictly stationary, reversible Markov chains, a connection between the ``geometric ergodicity'' condition and a certain ``spectral gap'' condition from functional analysis. 
It is well known that that ``key result'' can be formulated in the terminology of the dependence coefficients associated with the ``absolute regularity'' (``$\beta$-mixing'') condition and the ``$\rho$-mixing'' condition.
All that will be discussed in a little more detail after
Theorem \ref{th1.2} below.
In an effort to provide a relatively gentle exposition of 
that ``key result'' and some related  theory from the two aforementioned papers, 
the author \cite{ref-journal-Bradley2021} 
employed the dependence coefficients associated with the ``strong mixing'' (``$\alpha$-mixing'') condition as a means to try to organize that exposition into a somewhat more 
transparent form.

    That expository effort in \cite{ref-journal-Bradley2021} in turn leads to the question of what combinations of mixing rates, within some unavoidable narrow constraints, are possible for the 
$\rho$-mixing and strong mixing conditions together --- and perhaps including the absolute regularity condition as well --- for strictly stationary, reversible Markov chains.
This paper here is intended to give at least a partial answer to that question.  
The construction given in this paper --- a class of strictly stationary, reversible Markov chains --- will illustrate, within the unavoidable narrow constraints alluded to above, a pretty broad spectrum of possible combinations of such mixing rates.  As the examples will also involve a state space that is countable (rather than, say, the whole real number line), the arguments involved in checking the various properties of the examples here will be pretty elementary.

     To set up further discussion, Setting \ref{st1.1} below
will lay out some key terminology and background facts.
 
\medskip              
\begin{setting}\label{st1.1}   
First a few preliminary formalities:  
Throughout this paper, \hfil\break
$\R$ denotes the set of all real numbers, \hfil\break
$\Z$ denotes the set of all integers, and \hfil\break
$\N$ denotes the set of all positive integers.

   In this paper, a set $S$ is said to be ``countable'' if $S$ is either finite or countably infinite. 

   The usual notations (such as $\R^\Z$) will be used for Cartesian products of sets.
   
   The notation $\log x$ refers to the natural logarithm.

   If $(a_n, n \in \N)$ and $(b_n, n \in \N)$ are each a 
sequence of positive numbers, then the notation 
$a_n \asymp b_n$ as $n \to \infty$ will mean that 
$a_n = O(b_n)$ and $b_n = O(a_n)$ as $n \to \infty$.
   
    If $(a_n, n \in \N)$ is a sequence of nonnegative numbers, then the notation
``$a_n \to 0$ at least exponentially fast as $n \to \infty$'' means that there exists a positive number $c$ 
such that $a_n = O(e^{-cn})$ as $n \to \infty$.
\medskip    

     {\bf Part 1}\ \  ({\it The probability space}).\ \ The setting for the probability theory in 
this paper is a probability space $(\Omega, \cF, P)$, rich enough to accommodate all 
random variables declared.  
All random variables in this paper are defined on that probability space.

   For simplicity, {\it all random variables in this paper are real-valued\/}.  Some may be discrete, including those in the examples constructed in this paper.
(It is to be tacitly understood that some of the sources cited in this paper dealt with random variables taking their values in more general spaces than the real numbers; however, for 
the material in this paper, there is no need to deal with that greater generality.) 

     For a given family $(Y_i, i \in I)$ of (real-valued) random variables (where $I$ is a nonempty index set), the $\sigma$-field ($\subset \cF$, on $\Omega$) generated by this family will be denoted $\sigma(Y_i, i \in I)$. 
\medskip

     {\bf Part 2}\ \ ({\it Three measures of dependence\/}).\ \ Suppose $\cA$ and $\cB$ are any 
two $\sigma$-fields $\subset \cF$.
Define the following three measures of dependence:
First,
\begin{equation}\label{eq1.11}
\alpha(\cA,\cB) := \sup_{A\in \cA, B\in \cB} |P(A\cap B)-P(A) P(B)|.
\end{equation}
Next, 
\begin{equation}\label{eq1.12}
\beta(\cA,\cB):= \sup \frac{1}{2} \sum^I_{i=1}\sum^J_{j=1} |P(A_i\cap B_j) - P(A_i)P(B_j)|
\end{equation}
where the supremum is taken over all pairs of finite partitions 
$\{A_1,A_2,\dots,A_I\}$ and $\{B_1,B_2,\dots,B_J\}$ of 
$\Omega$ such that $A_i\in \cA$ for each $i$ and $B_j\in \cB$ for each~$j$.
[The factor of $1/2$ in (\ref{eq1.12}) is not of special significance, but has become customary in order to make 
certain inequalities ``a little nicer''.]\ \ 
Finally, define the ``maximal correlation coefficient''
\begin{equation}
\label{eq1.13}
\rho(\cA,\cB) := \sup | \corr(Y,Z)|
\end{equation}
where the supremum is taken over all pairs of square-integrable random variables $Y$ and $Z$ such that $Y$ is $\cA$-measurable and $Z$ is $\cB$-measurable.
\medskip 

The following inequalities are elementary and well known:
\begin{align}
\label{eq1.14}
0 &\le 2\alpha(\cA,\cB) \le \beta(\cA,\cB) \le 1; \quad {\rm and} \\
\label{eq1.15}
0 &\le 4\alpha(\cA,\cB) \le \rho(\cA,\cB) \le 1. 
\end{align}
(See e.g.\  [\cite{ref-journal-Bradley2007}, v.1, Proposition 3.11].)\ \  The quantities
$\alpha(\cA,\cB)$, $\beta(\cA,\cB)$, and $\rho(\cA,\cB)$ are all equal to $0$ if 
the $\sigma$-fields $\cA$ and $\cB$ are independent, and are all positive otherwise.
\medskip

     {\bf Part 3}\ \  ({\it Three mixing conditions\/}).\ \  Suppose $X:= (X_k$, $k\in \Z)$ is a (not necessarily Markovian) strictly stationary sequence of (real-valued) random variables.
For each integer $j$, define the notations $\cF^j_{-\infty} :=\sigma(X_k, k\le j)$ and 
$\cF^\infty_j:= \sigma(X_k$, $k\ge j)$.
\medskip

     For each positive integer $n$, define the following three ``dependence coefficients'':
\begin{align}
\label{eq1.16}
\alpha(n) &= \alpha_X(n) := \alpha(\cF^0_{-\infty}, \cF^\infty_n); \\
\label{eq1.17}
\beta(n) &= \beta_X(n) := \beta(\cF^0_{-\infty}, \cF^\infty_n);\,\,\mbox{and}\\
\label{eq1.18}
\rho(n) &= \rho_X(n) := \rho(\cF^0_{-\infty}, \cF^\infty_n).
\end{align}
For each positive integer $n$, one has by strict stationarity that 
$\alpha(n) = \alpha(\cF^j_{-\infty},\cF^\infty_{j+n})$ for every integer 
$j$, and the analogous comment holds for $\beta(n)$ and 
for $\rho(n)$ as well.

     Also, one (trivially) has that $\alpha(1) \ge \alpha(2)\ge \alpha(3)\ge \dots\,$; and the analogous comment holds for the numbers $\beta(n)$ and for the numbers $\rho(n)$.
\medskip

      The strictly stationary sequence $X$ is said to satisfy

\noindent ``strong mixing'' (or ``$\alpha$-mixing'') if $\alpha_X(n) \to 0$ as $n\to\infty$;

\noindent ``absolute regularity'' (or ``$\beta$-mixing'') if $\beta_X(n) \to 0$ as $n\to\infty$;

\noindent ``$\rho$-mixing'' if $\rho_X(n) \to 0$ as $n\to \infty$.
\medskip

The strong mixing ($\alpha$-mixing) condition is due to Rosenblatt \cite{ref-journal-Rosenblatt1956}.
The absolute regularity ($\beta$-mixing) condition was first studied by Volkonskii and Rozanov 
\cite{ref-journal-VolkRoz}, and was attributed there to Kolmogorov.
 The $\rho$-mixing condition is due to Kolmogorov and Rozanov \cite{ref-journal-KolRoz}.
(The ``maximal correlation coefficient'' $\rho(\cA,\cB)$ itself, for 
 $\sigma$-fields $\cA$ and $\cB$, was first studied earlier by Hirschfeld \cite{ref-journal-Hirschfeld} in a statistical context that had no particular connection with stochastic processes.)
\medskip 

    By (\ref{eq1.14}) and (\ref{eq1.15}), one has that for each positive integer $n$,
\begin{align}
\label{eq1.19}
0&\le 2\alpha(n) \le \beta(n) \le 1,\,\,\mbox{and}\\
\label{eq1.110}
0 &\le 4\alpha(n) \le \rho(n) \le 1.
\end{align}
By (\ref{eq1.19}), absolute regularity ($\beta$-mixing) implies strong mixing 
($\alpha$-mixing); and 
by (\ref{eq1.110}), $\rho$-mixing implies strong mixing ($\alpha$-mixing).   
[Part 5(E) below, and Example 1 of Remark \ref{rm1.4} 
later on, will together illustrate the fact that of the absolute regularity and 
$\rho$-mixing conditions, neither implies the other.]  
\medskip 

     {\bf Part 4}\ \ ({\it Strictly stationary Markov chains\/}).\ \ Now suppose that 
$X:=(X_k, k\in \Z)$ is a strictly stationary {\em Markov chain} (with the random variables $X_k$, $k \in \Z$ being real-valued, possibly discrete).
(No assumption of ``reversibility'' yet.)
\medskip

      As a well known consequence of the Markov property, for each positive integer $n$, eqs.\ (\ref{eq1.16})-(\ref{eq1.18}) hold in the following augmented forms for the given (strictly stationary) Markov chain $X$:
\begin{align}
\label{eq1.111}
\alpha(n) &= \alpha(\cF^0_{-\infty}, \cF^\infty_n) = \alpha\bigl(\sigma(X_0),\sigma(X_n)\bigl);\\
\label{eq1.112}
\beta(n) &= \beta(\cF^0_{-\infty},\cF^\infty_n) = \beta\bigl(\sigma(X_0), \sigma(X_n)\bigl);\\
\label{eq1.113}
\rho(n) &=\rho(\cF^0_{-\infty}, \cF^\infty_n) = \rho\bigl(\sigma(X_0), \sigma(X_n)\bigl).
\end{align}
(See e.g.\  [\cite{ref-journal-Bradley2007}, v.1, Theorem 7.3].)
\medskip

      By strict stationarity and (\ref{eq1.111})-(\ref{eq1.113}), one has that for any integer $j$ and any positive integer $n$, the (strictly stationary) Markov chain $X$ satisfies \hfil\break 
(i)\ $\alpha(n) = \alpha(\sigma(X_j),\sigma(X_{j+n}))$, \hfil\break
(ii)\ $\beta(n) = \beta(\sigma(X_j), \sigma(X_{j+n}))$, and \hfil\break 
(iii) $\rho(n) = \rho(\sigma(X_j), \sigma(X_{j+n}))$. \hfil\break
 
      Here are some special facts involving the 
dependence coefficients $\rho(n)$, $n \in \N$. 
As an elementary consequence of (\ref{eq1.13}), 
for any two $\sigma$-fields $\cA$ and $\cB$, 
$\rho(\cA, \cB) = \sup \|E(Y|\cB)\|_2/\|Y\|_2$ where the supremum is taken over all square-integrable, $\cA$-measurable
random variables $Y$ with mean 0.  (When necessary, interpret 0/0 := 0.)\ \ As a well known application of that fact
(together with the equality in (iii) in the preceding paragraph),
for the given strictly stationary Markov chain
$X:=(X_k, k\in \Z)$, one has that 
for any pair of positive integers $m$ and $n$, 
\begin{equation}
\label{eq1.114}
\rho(m+n) \le \rho(m) \cdot \rho(n).
\end{equation}
In particular, for any two positive integers $m$ and $n$, $\rho(m(n+1)) \le \rho(mn) \cdot \rho(m)$. Hence by induction, for every positive integer $m$, one has that the given (strictly stationary) Markov chain $X$ satisfies
\begin{equation}
\label{eq1.115}
\rho(mn) \le [\rho(m)]^n\,\,\mbox{for every $n\in\N$.}
\end{equation}
In particular (take $m=1$),
\begin{equation}\label{eq1.116}
\rho(n) \le [\rho(1)]^n\,\,\mbox{for every $n\in \N$.}
\end{equation}
By (\ref{eq1.115}), for the given (strictly stationary) Markov chain $X:= (X_k, k\in \Z)$, the following three conditions are equivalent: \hfil\break
(i) there exists $m\ge 1$ such that $\rho(m) <1$; \hfil\break
(ii) $X$ is $\rho$-mixing; \hfil\break
(iii) $\rho(n) \to 0$ at least exponentially fast as $n\to \infty$.
\medskip

   {\bf Part 5}\ \ ({\it Strictly stationary, reversible Markov chains\/}).  A given strictly stationary Markov chain $X := (X_k, k \in \Z)$ is said to be ``reversible'' if the distribution
(on $\R^{\Z}$) of the ``time-reversed'' sequence $(X_{-k}, k \in \Z)$ is identical to that of the sequence $X$ itself.
\medskip

   (A)\ \ By a well known argument, a given strictly stationary Markov chain
$X := (X_k, k \in \Z)$ is reversible if and only if the distribution (on $\R^2$) of the random vector $(X_1, X_0)$ is identical to that of the random vector $(X_0, X_1)$.
\medskip

   (B)\ \ If a given strictly stationary Markov chain 
$X := (X_k, k \in \Z)$ is reversible, then 
\begin{equation}\label{eq1.117}
\rho(n) = [\rho(1)]^n\,\,\mbox{for every $n\in \N$.}
\end{equation} 
That was shown by 
Longla [\cite{ref-journal-Longla14}, Lemma 2.1] in the context of strictly stationary, reversible Markov chains involving certain types of copulas; his argument extends beyond that context.
It can also be seen as an application of a certain theorem in functional analysis involving self-adjoint bounded linear operators from a Hilbert space to itself.
(This latter fact, explained in more detail in 
\cite{ref-journal-Bradley2021}, was pointed out and 
contributed there by an anonymous referee of that paper.
A presentation of the proof in probabilistic terminology was also included in the exposition in that paper.)\ \ 
Compare (\ref{eq1.117}) to (\ref{eq1.116})
(where reversibility was not assumed). 
\medskip

   (C)\ \    The paper \cite{ref-journal-Bradley2015} 
constructed some strictly stationary, countable-state, reversible Markov chains that are $\rho$-mixing but fail to satisfy  $\rho^*$-mixing --- the ``interlaced'' variant of $\rho$-mixing with the two index sets not restricted to ``past'' and ``future''.
In those examples, 
$\rho(\sigma(X_0), \sigma(X_{-n}, X_n)) = 1$
for every $n \in \N$.
\medskip 

  (D)\ \ As is well known, any strictly stationary 
{\it two-state\/} Markov chain
$X := (X_k, k \in \Z)$ is reversible.
\smallskip

   The explanation is simple.   
Let us label the two states as 0 and 1.
Then by strict stationarity,
\begin{align*}
P(\{X_0 = 0\} \cap \{X_1 = 1\})\ 
&=\ P(X_0 = 0) - P(\{X_0 = 0\} \cap \{X_1 = 0\})\\
&=\ P(X_1 = 0) - P(\{X_0 = 0\} \cap \{X_1 = 0)\}\
=\ P(\{X_0 = 1\} \cap \{X_1 = 0\}). 
\end{align*}
That is, $P((X_0, X_1) = (0,1)) = P((X_1, X_0) = (0,1))$.
Now it trivially follows that
$P((X_0, X_1) = (i,j)) = P((X_1, X_0) = (i,j))$
for all ordered pairs $(i,j) \in \{0,1\}^2$.
Hence by (A) above, one has that (D) holds.
\medskip

   (E)\ \ There exist strictly stationary, reversible Markov chains $X := (X_k, k \in \Z)$ that satisfy (i) $\rho$-mixing, and 
(ii) $\beta(n) = 1$ for all $n \in \N$. 
Condition (ii) is (for strictly stationary Markov chains) equivalent to the condition that for every positive integer $n$, the $n$-step transition distributions are almost surely totally singular with respect to the (invariant) marginal distribution.  
Such examples cannot be countable-state.
\smallskip

     Using ``random rotations'', Rosenblatt [\cite{ref-journal-Rosenblatt1971}, pp.\ 214-215] constructed some strictly stationary Markov chains that satisfy conditions (i) and (ii) above.
It seems clear that some of those examples (the ones for which the underlying ``random rotations'' satisfy an appropriate symmetry) are reversible.
As a slight variant of those examples of Rosenblatt,
the author 
[\cite{ref-journal-Bradley2007}, Examples 7.16 and 7.17] 
constructed examples of strictly stationary, reversible Markov chains that satisfy (i) and (ii) above, and showed that 
those examples in fact satisfy the $\rho^*$-mixing condition 
alluded to in (C) above.
Reversibility was not mentioned there (for those examples 
in \cite{ref-journal-Bradley2007}), but is easily verified as a consequence of the reversibility of the ``building blocks'' for those examples --- strictly stationary, two-state (hence reversible) Markov chains (see (D) above). 
\end{setting}     

   The following known theorem will help focus the further discussion.  \medskip

\begin{theorem}
\label{th1.2}
\textit{Suppose $X := (X_k,\ k\in \Z)$ is a strictly stationary, reversible Markov chain, and
$r$ is a real number such that $0 < r < 1$.
Then the following four conditions are equivalent: \hfil\break 
(i) $\rho_X(1) \leq r$. \zhfb
(ii) For all $n \in \N$, $\rho_X(n) \leq r^n$. \zhfb
(iii) For all $n \in \N$, $\alpha_X(n) \leq r^n$. \zhfb
(iv) $\alpha_X(n) =  O(r^n)$ as $n \to \infty$. }
\end{theorem}

     Except for one superficial simplification here, this theorem is [\cite{ref-journal-Bradley2021},Corollary 5.7].  
Obviously (i) $\Leftrightarrow$ (ii) by Part 5(B) of Setting \ref{st1.1}; and (see (\ref{eq1.110})) trivially 
(ii) $\Rightarrow$ (iii) $\Rightarrow$ (iv).
The proof of the ``remaining implication'' (iv) $\Rightarrow$ (i) is the argument that was given for 
[\cite{ref-journal-Bradley2021}, Lemma 5.5] --- an argument that was just a slight variant of calculations in the papers of Roberts and Rosenthal \cite{ref-journal-RR} 
and Roberts and Tweedie \cite{ref-journal-RT}.

   The theory developed by Roberts, Rosenthal, and Tweedie in those two papers \cite{ref-journal-RR} and \cite{ref-journal-RT} 
contains (in greater generality, and together with other things) the ``key result'' alluded to earlier.
Without formally defining all of the necessary terminology, 
one can at least state that ``key result'' informally as follows:
{\it If a given strictly stationary Markov chain is reversible and satisfies a certain ``irreducibility'' condition (equivalent to Harris recurrence), then the ``geometric ergodicity'' condition holds if and only if a certain ``$L^2$ spectral gap'' condition from functional analysis holds.}\ \
As noted earlier, it is well known that that ``key result'' and other related results can be transcribed into the terminology of dependence coefficients. 
See for example 
\cite{ref-journal-KM}, 
\cite{ref-journal-Longla14}, and 
\cite{ref-journal-LP},
or the detailed exposition in \cite{ref-journal-Bradley2021}. 
The papers of 
Nummelin and Tweedie \cite{ref-journal-NTweed} and 
Nummelin and Tuominen \cite{ref-journal-NTuo} 
together had (among other things) already shown (in different terminology) that for a given strictly stationary Markov chain (reversible or not), the geometric ergodicity condition is equivalent to absolute regularity with 
$\beta(n) \to 0$ at least exponentially fast.
(For more on geometric ergodicity, see e.g.\ 
\cite{ref-journal-MT} or
[\cite{ref-journal-Bradley2007}, v.2, Chapter 21].)\ \ 
It is elementary and long well known that the 
``$L^2$ spectral gap'' condition alluded to above is, for a given strictly stationary Markov chain (reversible or not) equivalent to
the condition $\rho(1) < 1$. 
The ``key result'' of Roberts, Rosenthal, and Tweedie from above, together with a ``comparison of rates of convergence'' calculation of theirs connected with it, can
together be formulated in a way that is ``almost'' (but ``not quite'') an analog of Theorem \ref{th1.2} with $\beta(n)$ in place of $\alpha(n)$.

     All that is reviewed in a bit of detail in 
\cite{ref-journal-Bradley2021}, where in essence Theorem \ref{th1.2} was developed mainly in order to assist an expository effort to provide a gentle introduction to the ``key result'' of Roberts, Rosenthal, and Tweedie.
In contrast to that ``key result'', Theorem \ref{th1.2} does not require any assumption of ``irreducibility'' (or Harris recurrence); it applies just as well to the peculiar $\rho$-mixing examples alluded to in Part 5(E) of Setting \ref{st1.1}, where the property $\beta(n) = 1$ for all $n \in \N$ actually prevents such ``irreducibility'' or Harris recurrence (and prevents the state space from being countable).

     In this paper here, with Theorem \ref{th1.2} as the background starting point, we shall investigate the question of what combinations of behavior of the dependence coefficients 
$\alpha(n)$ and $\rho(n)$ are possible for strictly stationary, reversible Markov chains.

     In addition to the sources cited above in connection with the
``key result'' of Roberts, Rosenthal, and Tweedie, the main sources of inspiration for this paper here are some earlier papers
that (among other things) established some  
``slower than exponential'' mixing rates for the strong mixing condition for certain classes of strictly stationary (not necessarily reversible) Markov chains with a ``renewal structure''.    
To illustrate certain limitations in connection 
with central limit theorems for strictly stationary 
(not necessarily Markovian), strongly mixing sequences,  
Davydov \cite{ref-journal-Davydov},
Tikhomirov \cite{ref-journal-Tikh},
and Doukhan, Massart, and Rio \cite{ref-journal-DMR}
each constructed strictly stationary Markov chains for which,
for some appropriate positive number $c$, 
$\alpha(n) \asymp \beta(n) \asymp n^{-c}$ as $n \to \infty$.
In those papers, the main focus was on $\alpha(n)$. 
For the examples in \cite{ref-journal-Davydov} and
\cite{ref-journal-Tikh}, the inclusion of $\beta(n)$ was implicitly established by the calculations of 
Davydov [\cite{ref-journal-Davydov}, pp.\ 327-328]; 
and in the paper
\cite{ref-journal-DMR} the inclusion of $\beta(n)$ was
established explicitly.  
Kesten and O'Brien [\cite{ref-journal-KO}, pp.\ 412-414] constructed a broad class of examples of strictly stationary, strongly mixing Markov chains with a quite arbitrary 
``slower than exponential'' mixing rate.
The Markov chains there (as well as those in
\cite{ref-journal-Davydov} and
\cite{ref-journal-Tikh}) were countable-state.
The Markov chains constructed in \cite{ref-journal-DMR} 
are reversible (and have a ``continuum'' state space, namely  the unit interval $[0,1]$).  
  
   This paper here is focused partly on extending the observations of Kesten and O'Brien \cite{ref-journal-KO} cited above (displaying a broad spectrum of possible ``sub-exponential'' mixing rates for $\alpha(n)$ for strictly 
stationary Markov chains)
to strictly stationary, Markov chains that are reversible. 
This paper is also intended to give, for strictly stationary 
Markov chains that are reversible and $\rho$-mixing, similar results on a ``broad spectrum of pairs of mixing rates'' (to the narrow extent permitted by Part 5(B) of Setting \ref{st1.1} and by Theorem 1.2) involving the dependence coefficients $\rho(n)$  and $\alpha(n)$ together.
The Markov chains constructed in this paper will be 
countable-state; and they will include absolute regularity, with
$\alpha(n) \asymp \beta(n)$ as $n \to \infty$.
  
   By (\ref{eq1.117}) and the second sentence after (\ref{eq1.15}), for a given strictly stationary, reversible Markov chain, the dependence coefficients $\alpha(n)$, $\beta(n)$, 
and $\rho(n)$, $n \in \N$, will {\it all\/} be positive -- except in the case of a sequence of independent, identically distributed random variables (where those dependence coefficients are all zero).
In the ``dependent'' cases, i.e.\ where $0 < \rho(1) \leq 1$ (the value 1 is included here), the ratio $\alpha(n)/\rho(n)$ is bounded above by $1/4$, by (\ref{eq1.110}).
The focus in this paper will be on cases where that ratio converges to 0 as $n \to \infty$.

   Let it first be noted in passing that 
[see Part 5(D) of Setting \ref{st1.1}] for some well known strictly stationary, two-state, (hence) reversible Markov chains, the strong mixing, $\rho$-mixing, and absolute regularity conditions hold with the same (exponential) mixing rate 
modulo a constant factor.
In Section \ref{sc2}, that will be pointed out in concrete simple detail in Lemma \ref{lm2.5} and its subsequent Remark.

   Now (in the ``dependent'' case) as a consequence of 
Theorem \ref{th1.2} (with no changes), one has that whether 
$0 < \rho_X(1)  < 1$ or $\rho_X(1) = 1$, if the ratio
$\alpha(n)/\rho(n)$ converges to 0, it cannot do so at an (at least) exponential rate.
The class of examples that will be presented in this paper, will show that that ratio can converge to 0 at a quite arbitrary ``slower than exponential'' rate, at least within a mild ``log convexity'' condition --- regardless of whether $\rho$-mixing holds or instead $\rho(n) = 1$ for all $n \in \N$
[see (\ref{eq1.117}) or again the last sentence 
(i.e.\ the equivalence) of Part 4 of Setting \ref{st1.1}]. 
   
    Here is our main result (stated here with some redundancy):

\begin{theorem}
\label{th1.3}
\textit{Suppose $r$ is a real number such that $0 < r \leq 1$.
Suppose $f: [0, \infty) \to (0, 1/2]$ is a continuous, strictly decreasing function with the following three properties: \zhfb 
(i)\ \ $f(x) \to 0$ as $x \to \infty$; \zhfb
(ii)\ \ for every $u \in (0,1)$, $u^x = o(f(x))$ as
$x \to \infty$; and \zhfb
(iii)\ \ the mapping $x \mapsto \log f(x)$, for 
$x \in [0,\infty)$, is convex on $[0, \infty)$. \zhfb
\indent  Then there exists (on the probability space 
$(\Omega, \cF, P)$) a strictly stationary,
countable-state Markov chain $X := (X_k,\ k\in \Z)$ 
such that $X$ is reversible and for every positive integer $n$,
\begin{equation}
\label{eq1.33}
\rho_X(n) = r^n \quad {\rm and} \quad 
(1/2)\, r^n f(n) \leq \alpha_X(n) \leq \beta_X(n) 
\leq 12\, r^n f(n). 
\end{equation}  
}
\end{theorem}

    Here the numbers $1/2$ (twice) and 12 are not of special significance, but are instead just an artifact of an endeavor
to make the statement and proof of Theorem 1.3 a ``little nicer''.  [Again recall (\ref{eq1.110}).]\ \    
Theorem \ref{th1.3} will be proved in Section \ref{sc3},
after some preliminary work is done in Section \ref{sc2}.
The rest of Section \ref{sc1} here will be devoted to
certain specific illustrations of Theorem \ref{th1.3}.

\medskip  
\begin{remark}
\label{rm1.4} Here we shall take a quick look at a couple of related specific applications of 
Theorem \ref{th1.3}.
In the first, with $r=1$, $\rho$-mixing fails to hold; in the 
second, with $0 < r < 1$, $\rho$-mixing holds.
\medskip 

   {\it Example 1.}\ \ This example is in spirit an adaptation,
to reversibility, of an illustration in the work of
Kesten and O'Brien \cite{ref-journal-KO} alluded to above.  
   
   Suppose that $0 < a< 1$ and $q > 0$, and also that $b$ and $c$ are each an arbitrary real number.
Define the function $\eta: (1,\infty) \to \R$ as follows:
\begin{equation*}
{\rm for\ all}\ x \in (1,\infty),\ \ 
\eta(x) := \bigl(\exp(-qx^a)\bigl)\, \cdot\, x^b 
\cdot (\log x)^c. 
\end{equation*}
Then $\eta(x) \to 0$ as $x \to \infty$.
Define the function $h: (1,\infty) \to \R$ as follows:
\begin{equation*}
{\rm for\ all}\ x \in (e,\infty),\ \ 
h(x) := \log(\eta(x)) = -qx^a + b (\log x) + c \log(\log x). 
\end{equation*}
Then by elementary calculations, the following statements hold: \zhfb
(i)\ $\eta(x) \to 0$ as $x \to \infty$; 
and for all $t > 0$, 
$\eta(t+x)/\eta(x) \to 1$ as $x \to \infty$; \zhfb 
(ii)\ $h'(x) \to 0$ as $x \to \infty$; and \zhfb
(iii)\ for all sufficiently large $x > 0$, one has that 
$h'(x) < 0$ and $h''(x) > 0$.
\smallskip  

     From all of the above observations on the functions 
$\eta$ and $h$, one has the following:
If $T > 0$ is chosen sufficiently large, the function 
$f: [0, \infty) \to \R$ defined by
$f(x) = \eta(T + x)$ for $x \in [0, \infty)$ will be well defined and will satisfy all of the hypotheses (for $f$) in Theorem \ref{th1.3}.
Consequently, by that theorem (with $r = 1$ there), 
there exists a strictly stationary, reversible Markov chain 
$X:= (X_k, k \in \Z)$ (with countable state space), such that
$\rho(n) = 1$ for all $n \in \N$, and
strong mixing and absolute regularity
both hold with mixing rates
\begin{equation*}
\alpha_X(n) \asymp \beta_X(n)  \asymp 
\bigl(\exp(-qn^a)\bigl)\, \cdot\, n^b \cdot (\log n)^c\ \ 
{\rm as}\ n \to \infty.  
\end{equation*} 

   {\it Example 2.}\ \ Now for a given $r$ such that $0 < r < 1$, and for $a$, $q$, $b$, and $c$ with the same restrictions as in Example 1 above, this time applying Theorem \ref{th1.3} with the given $r$ and
(again) $f(x) = \eta(T + x)$ where $\eta$ and $T$ are as above, one has that there exists a strictly stationary, countable-state, reversible Markov chain $X:= (X_k, k \in \Z)$ such that  
$\rho$-mixing holds with $\rho_X(n) = r^n$ for all $n \in \N$, and the strong mixing and absolute regularity conditions hold with mixing rates
 \begin{equation*}
\alpha_X(n) \asymp \beta_X(n)  \asymp r^n \cdot 
\bigl(\exp(-qn^a)\bigl)\, \cdot\, n^b \cdot (\log n)^c\ \ 
{\rm as}\ n \to \infty. 
\end{equation*} 
\end{remark}

\section{Preliminaries}\label{sc2}

\begin{notations}\label{nt2.1}
The construction (in Section \ref{sc3}) of the strictly stationary, countable-state, reversible Markov chain for Theorem \ref{th1.3} will involve as ``building blocks'' a countably infinite collection of strictly stationary, 2-state, reversible Markov chains
that are independent of each other.
(Again recall Part 5(D) of Setting \ref{st1.1}.)\ \ 
The arithmetic for that process apparently will be slightly less ``cluttered'' if in those 2-state
``building block'' Markov chains, the states are labeled 0 and 1 (instead of, say, 1 and 2).
Accordingly, the $2 \times 2$ matrices involved in that process (for example, the $n$-step
transition probability matrices) will be set up in the form $M := (m_{ij})_{i,j \in \{0,1\}}$, with the
top row being $[m_{00}, m_{01}]$ and the bottom row being $[m_{10}, m_{11}]$.
\medskip

   {\bf Part 1.}\ \ {\it Joint probability matrices\/}.\ \ 
First we shall spell out a class of $2 \times 2$ matrices that will be used for {\it joint\/} 
(not transition) probabilities.  

     For any $\zeps \in (0, 1/2]$ and any $\theta \in (0,1)$, define the $2 \times 2$ matrix
$\Lambda^{(\zeps, \theta)} := (\lambda^{(\zeps,\theta)}_{ij})_{i,j \in \{0,1\}}$ as follows:
\begin{align}
\label{eq2.11}
\lambda^{(\zeps,\theta)}_{00} &:= (1 - \zeps)^2 + (1 - \zeps) \zeps \theta;\\
\label{eq2.12}
\lambda^{(\zeps,\theta)}_{01} = \lambda^{(\zeps,\theta)}_{10} 
&:= (1 - \zeps) \zeps - (1 - \zeps) \zeps \theta;\ \ {\rm and}\\
\label{eq2.13}
\lambda^{(\zeps,\theta)}_{11} &:= \zeps^2 + (1 - \zeps) \zeps \theta.
\end{align}
Note that (under the stipulated conditions on $\zeps$ and $\theta$), these four entries are all positive, and their sum is 1.

We shall return to that class of matrices in Lemma \ref{lm2.2}\ below.
\medskip

     {\bf Part 2.}\ \  {\it Transition probability matrices\/}.\ \ The matrices below will play a key role in transition probabilities (of some strictly stationary, 2-state Markov chains).
     
     The $2 \times 2$ identity matrix  $(\delta_{ij})_{i,j \in \{0,1\}}$ (where
$\delta_{00} = \delta_{11} = 1$ and $\delta_{01} = \delta_{10} = 0$)  
will be denoted simply as $I_2$.

     For any $\zeps \in (0,1/2]$, let $A^{(\zeps)} :=  (a^{(\zeps)}_{ij})_{i,j \in \{0,1\}}$ be the
 $2 \times 2$ matrix in which each of the two rows is $[1 - \zeps, \zeps]$.

     For any $\zeps \in (0,1/2]$ and any $\theta \in (0,1)$, define the $2 \times 2$ matrix
$\bbP^{(\zeps,\theta)} := (p^{(\zeps,\theta)}_{ij})_{ i,j \in \{0,1\}}$ as follows:
\begin{equation}
\label{eq2.14}
\bbP^{(\zeps,\theta)} := \theta I_2 + (1 - \theta) A^{(\zeps)} .      
\end{equation}
By simple arithmetic, for any $\zeps \in (0,1/2]$ and any 
$\theta \in (0,1)$, the entries of the
$2 \times 2$ matrix $\bbP^{(\zeps,\theta)}$ are as follows:
\begin{align}
\label{eq2.15}
p^{(\zeps, \theta)}_{00} &= (1 - \zeps) + \zeps \theta\ \ \ {\rm and}\ \ \ 
p^{(\zeps, \theta)}_{01} = \zeps - \zeps \theta;\\
\label{eq2.16} 
p^{(\zeps, \theta)}_{10} &= (1 - \zeps) - (1 - \zeps) \theta\ \ \ {\rm and}\ \ \ 
p^{(\zeps, \theta)}_{11} = \zeps + (1 - \zeps) \theta.
\end{align}
Note that in each of the two rows of $\bbP^{(\zeps,\theta)}$, the two entries are positive and their sum is 1.

   By trivial matrix multiplication, $(A^{(\zeps)})^2 = A^{(\zeps)}$ for any $\zeps \in (0,1/2]$.   
Hence by simple matrix multiplication, for any $\zeps \in (0,1/2]$ and any pair of numbers
$\theta, \tau \in (0,1)$, one has that
\begin{equation}
\label{eq2.17}
\bbP^{(\zeps, \theta)} \bbP^{(\zeps, \tau)} = \bbP^{(\zeps, \theta \tau)}.
\end{equation} 

   We shall return to those notations in Lemma \ref{lm2.3} below.
\bigskip
   
    {\bf Part 3.}\ \ {\it Alternative notations\/}.\ \ In the case where the number $\theta \in (0,1)$ itself involves ``smaller print'', the matrices $\Lambda^{(\zeps, \theta)}$ and
$\bbP^{(\zeps, \theta)}$ may be written as
$\Lambda(\zeps, \theta)$ and $\bbP(\zeps,\theta)$ respectively for typographical convenience.
\end{notations}

\begin{lemma}
\label{lm2.2}
Suppose $\zeps \in (0, 1/2]$ and $\theta \in (0,1)$.
Suppose $Y$ and $Z$ are $\{0,1\}$-valued random variables such that the (joint)
probability function of the random vector $(Y,Z)$ is the matrix 
$\Lambda^{(\zeps, \theta)}$
in (\ref{eq2.11})-(\ref{eq2.13}), that is, 
\begin{equation}
\label{eq2.21} 
{\rm for\ all}\ (i,j) \in \{0,1\}^2,\ \ \ P\bigl((Y,Z) = (i,j)\bigl) = \lambda^{(\zeps, \theta)}_{ij}\ . 
\end{equation}
Then the following statements hold:

     (I)\ \ The (joint) probability function of the random vector 
$(Z,Y)$ is the same as that of $(Y,Z)$.

     (II)\ \ $P(Y=0) = P(Z=0) = 1 - \zeps$\ \  and\ \  $P(Y = 1) = P(Z=1) = \zeps$.
     
     (III)\ \ $\corr(Y,Z) = \theta$.
     
     (IV)\ \ $\rho(\sigma(Y), \sigma(Z)) = \theta$.
     
     (V)\ \ $\alpha(\sigma(Y), \sigma(Z)) = (1 - \zeps) \zeps \theta \geq \zeps \theta /2$.  
  
     (VI)\ \ $\beta(\sigma(Y), \sigma(Z)) = 2(1 - \zeps) \zeps \theta \leq 2 \zeps \theta$. 
      
\end{lemma}

   {\bf Proof.}\ \ To verify statement (I), it in fact suffices to note that 
$P((Y,Z) = (0,1)) = P((Y,Z) = (1,0))$ by (\ref{eq2.21}) and (\ref{eq2.12}).

   Statement (II) holds by (\ref{eq2.21}) and (\ref{eq2.11})-(\ref{eq2.13}) and simple arithmetic.  

   To verify statement (III), use statement (II), 
eq.\ (\ref{eq2.21}), eqs.\ (\ref{eq2.11})-(\ref{eq2.13}), and simple arithmetic to show that 
$\var(Y) = \var(Z) = (1- \zeps)\zeps$ and
$\cov(Y,Z) = (1 - \zeps)\zeps \theta$.  
Statement (III) then follows.

   {\it Proof of (IV)\/.}\ \ Suppose $V$ and $W$ are any two random variables such that
$V$ (resp.\ $W$) is measurable with respect to $\sigma(Y)$ (resp.\ $\sigma(Z))$.
Since $Y$ and $Z$ each take just two values (0, and 1), it is easy to show that
there exist real numbers $a,\, b, \, c, \, d$ such that $V = aY + b$ and $W = cZ + d$.
One then has by elementary properties of correlation that
$\corr(V,W) = \corr(Y,Z)$ (resp.\ 0, resp.\ $-\corr(Y,Z)$) if $ac > 0$   
(resp.\ $ac=0$, resp.\ $ac < 0$).  
Then (IV) follows by (III) and eq.\ (\ref{eq1.13}) (and the fact that $\theta > 0$).

     {\it Proof of (V)-(VI)\/.}\ \ If $A$ and $B$ are events, and either one of them has probability 0 or 1, then trivially $P(A \cap B) - P(A)P(B) = 0$.  
Hence in using (\ref{eq1.11})-(\ref{eq1.12}) to calculate the dependence coefficients
$\alpha(\sigma(Y), \sigma(Z))$ and $\beta(\sigma(Y), \sigma(Z))$, one can restrict to pairs of events (say $A$ and $B$) whose probabilities are each strictly between 0 and 1.
By statement (II), eq.\ (\ref{eq2.21}), 
eqs.\ (\ref{eq2.11})-(\ref{eq2.13}), and simple arithmetic, one has that for every ordered pair $(i,j) \in \{0,1\}^2$,
\begin{equation}
\label{eq2.22}
\bigl|P(\{Y=i\} \cap \{Z=j\}) - P(Y=i)P(Z=j)\bigl|\, 
= (1 -\zeps) \zeps \theta.
\end{equation} 
From  (\ref{eq1.11})-(\ref{eq1.12}), one has that statements (V) amd (VI) both follow immediately from (\ref{eq2.22}) and the fact that $1/2 \leq 1 - \zeps < 1$.
That completes the proof of Lemma \ref{lm2.2}.
\medskip

\begin{lemma}
\label{lm2.3}
Suppose $\zeps \in (0, 1/2]$ and $\theta \in (0,1)$.
Suppose $V$ and $W$ are $\{0,1\}$-valued random variables such that
\begin{equation}
\label{eq2.31}
P(V = 0) = 1 - \zeps\ \ \ {\rm and}\ \ \ P(V=1) = \zeps. 
\end{equation}
and (see (\ref{eq2.15})-(\ref{eq2.16})) for each ordered pair $(i,j) \in \{0,1\}$,
\begin{equation}
\label{eq2.32}
P(W=j | V=i) = p^{(\zeps, \theta)}_{ij}
\end{equation}
(the $ij$-entry of the matrix $\bbP^{(\zeps, \theta)}$). 

   Then the (joint) probability function of the random vector $(V,W)$ is the 
matrix $\Lambda^{(\zeps, \theta)}$ (see (\ref{eq2.11})-(\ref{eq2.13})); and also in particular, 
\begin{equation}
\label{eq2.33}
P(W = 0) = 1 - \zeps\ \ \ {\rm and}\ \ \ P(W=1) = \zeps. 
\end{equation}      
\end{lemma}

   {\bf Proof.}\ \ For each ordered pair $(i,j) \in \{0,1\}^2$, one uses 
(\ref{eq2.31}) and (\ref{eq2.32}) together 
to show that (see (\ref{eq2.11})-(\ref{eq2.13})) 
$P((V,W) = (i,j)) = \lambda^{(\zeps, \theta)}_{ij}$.
After that, (\ref{eq2.33}) follows by trivial arithmetic (or from Lemma \ref{lm2.2}(II)).
Thus Lemma  \ref{lm2.3}\ holds.
\medskip

\begin{definition}
\label{df2.4}
Suppose $\zeps \in (0, 1/2]$ and $\theta \in (0,1)$.
A given random sequence $X := (X_k, k \in \Z)$ is is said to satisfy ``Condition $\cS(\zeps, \theta)$'' if $X$ is  a strictly stationary Markov chain with the following three properties:

   (1)\ \ The state space of $X$ is $\{0,1\}$.
   
   (2)\ \ The (invariant) marginal distribution is given by
\begin{equation}
\label{eq2.41}
P(X_0 = 0) = 1 - \zeps\ \ \ {\rm and}\ \ \ P(X_0=1) = \zeps.    
\end{equation}

   (3)\ \ The one-step transition probability matrix for $X$ is $\bbP^{(\zeps, \theta)}$
[see (\ref{eq2.14}) and (\ref{eq2.15})-(\ref{eq2.16})].         
\end{definition} 

   Condition $\cS(\zeps, \theta)$ is well defined.
By Lemma \ref{lm2.3}, the one-step transition probability matrix $\bbP^{(\zeps, \theta)}$
in (3) is compatible with the (invariant) marginal distribution in (\ref{eq2.41}).
\medskip

\begin{lemma}
\label{lm2.5}
Suppose $\zeps \in (0, 1/2]$ and $r \in (0,1)$.
Suppose $X := (X_k, k \in \Z)$ is a strictly stationary Markov chain that satisfies
Condition $\cS(\zeps,r)$. 
Then the following statements hold:

     (I)\ \ The Markov chain $X$ is reversible.
     
     (II)\ \ For each positive integer $n$, the $n$-step transition probability matrix for $X$ is
$\bbP(\zeps, r^n)$ (see eq.\ (\ref{eq2.14}), eqs.\ (\ref{eq2.15})-(\ref{eq2.16}), and Part 3 of  
Notations \ref{nt2.1}).

     (III)\ \ For each positive integer $n$, the following holds: 
\begin{align}
\label{eq2.51}
\rho_X(n) &= r^n;\\ 
\label{eq2.52}
\alpha_X(n) &= (1 - \zeps) \zeps r^n \geq \zeps r^n/2;\ \  
{\rm and}\\ 
\label{eq2.53}
\beta_X(n) &= 2(1 - \zeps) \zeps r^n \leq 2 \zeps r^n. 
\end{align}      
\end{lemma}

     {\bf Remark.}\ \  
    By (III), for any choice of parameters 
$\zeps \in (0, 1/2]$ and $r \in (0,1)$,
the (strictly stationary, two-state, reversible) Markov chain $X$ here in Lemma \ref{lm2.5} has the property that 
the dependence coefficients $\rho(n)$, $\alpha(n)$, and
$\beta(n)$ all converge to 0 with the same (exponential)
mixing rate, modulo constant factors. 
\medskip   

   {\bf Proof of Lemma \ref{lm2.5}.}\ \  
By Definition \ref{df2.4}, the (strictly stationary) Markov chain $X$ has state space $\{0,1\}$, with the marginal distribution given by (\ref{eq2.41}).
Statement (I) holds by Part 5(D) of Setting \ref{st1.1}.

   {\it Proof of (II).}\ \ By condition (3) in Definition \ref{df2.4}, one has that for each positive integer $n$, the $n$-step transition probability matrix for the Markov chain $X$ is
$(\bbP^{(\zeps,r)})^n$ (the $n^{\rm th}$ power with matrix multiplication).    
By (\ref{eq2.17}) and induction, one has that (see 
Part 3 of Notations \ref{nt2.1})    
$(\bbP^{(\zeps,r)})^n = \bbP(\zeps, r^n)$ for every positive integer $n$.
Thus (II) holds.

   {\it Proof of (III).}\ \ By statement (II) and Lemma \ref{lm2.3}, one has that for each
positive integer $n$, the (joint) probability function of the random vector $(X_0, X_n)$
is (again see Part 3 of Notations \ref{nt2.1}) the matrix 
$\Lambda(\zeps,r^n)$.
By (\ref{eq1.113}) and Lemma \ref{lm2.2}(IV), for each positive integer $n$,
$\rho_X(n) = \rho(\sigma(X_0), \sigma(X_n)) = r^n$. 
Thus (\ref{eq2.51}) holds. 
Eqs.\ (\ref{eq2.52})-(\ref{eq2.53}) follow similarly from (\ref{eq1.111})-(\ref{eq1.112}) and
Lemma \ref{lm2.2}(V)(VI).
Thus (III) holds.
That completes the proof of Lemma \ref{lm2.5}.
\medskip

     We shall return to Definition \ref{df2.4} and Lemma \ref{lm2.5} in Section \ref{sc3}.
Here in Section \ref{sc2}, we shall finish with a technical statement that will be employed in
Section \ref{sc3} as part of the process of assembling together a countably infinite collection of strictly stationary, two-state Markov chains (all independent of each other) from
Definition \ref{df2.4} in order to form the Markov chain for Theorem \ref{th1.3}.               
\medskip

\begin{lemma}
\label{lm2.6}
Suppose $I$ is a nonempty, countable (index) set, and 
$\cA_i,\, i \in I$ and $\cB_i,\, i \in I$ are $\sigma$-fields such that
the $\sigma$-fields $\cA_i \vee \cB_i,\, i \in I$ are independent.   
Then
\begin{align}
\label{eq2.61}
\beta \Bigl(\ \bigvee_{i \in I} \cA_i,\ \bigvee_{i \in I} \cB_i \Bigl)\
&\leq\ \sum_{i \in I} \beta(\cA_n, \cB_n);\ \ \ {\rm and}\\ 
\label{eq2.62}
\rho \Bigl(\ \bigvee_{i \in I} \cA_n,\ \bigvee_{i \in I} \cB_n \Bigl)\
&=\ \sup_{i \in I} \rho(\cA_n, \cB_n).
\end{align}      
\end{lemma}

   The latter equality\ (\ref{eq2.62}) is due to Cs\'aki and Fischer \cite{ref-journal-CF}.
   Both (\ref{eq2.61}) and (\ref{eq2.62}) (and under the same assumptions, the inequality
(\ref{eq2.61}) with $\beta$ replaced on both sides by $\alpha$, an inequality that will not  
be used here) can be found in
[\cite{ref-journal-Bradley2007}, v.1, Theorem 6.2].  
 [Actually, by elementary arguments, all results in that particular theorem, including the ones cited here, extend 
(pretty frivolously) to the case of an arbitrary (not necessarily countable) nonempty index set; that will not be needed here.]

\section{Proof of Theorem 1.3}\label{sc3}

     In this proof, in order to avoid some unnecessary clutter, the use of properties of the function $f$ in Theorem 1.3 (such as the ``log convexity'' assumption) will sometimes be at ``less than full strength'', resulting in some crude inequalities and leaving unsaid some natural ``extra'' or ``stronger'' observations that will not be needed in the argument.
     
     In this proof, just in an attempt to ``make the arithmetic a little nicer'', base-2 logarithms will be used, instead of natural logarithms.    

     The proof will be spelled out here in a series of (mostly) 
small ``steps'', identified by consecutive letters
(``Step A'', ``Claim B'', ``Step C'', and so on). 

   A large portion of this proof here will involve setting up some
elementary but tedious  ``scaffolding'' involving lines in connection with the ``log convexity'' assumption in the theorem.
Essentially the same scaffolding, with just minor differences, was used by the author \cite{ref-journal-Bradley1987}
to prove a quite different but somewhat related theorem (involving more mixing conditions, in a
non-Markovian, ``non-reversible'' context). 
There the scaffolding was spelled out tersely.
In the presentation of that theorem and its proof (again) 
in [\cite{ref-journal-Bradley2007}, v.3, Theorem 26.5],
that scaffolding was spelled out in detail, partly in the proof 
itself and partly in the Appendix at the end of that book.
As compared to the context there, the context here will involve some small but nontrivial differences in some information that is developed and used.   
For the reader's convenience, we shall spell out the scaffolding here in detail, in a self-contained form that is designed to 
cleanly fit our context here.  
\medskip

     {\bf Step A.}\ \  As in the statement of Theorem 1.3, suppose $0 < r \leq 1$, and 
$f:[0,\infty) \to (0, 1/2]$ is a continuous, strictly decreasing function satisfying assumptions (i), (ii), and (iii) in the statement
of that theorem. 

     Define the continuous, strictly decreasing (and negative) function
$g: [0,\infty) \to (-\infty, -1]$ as follows:
 \begin{equation}
\label{eq3A1}
     {\rm For\ each}\ x \in [0,\infty),\ \ \ 
g(x)\ :=\ \log_2 [\, r^x\, f(x)]\ =\ x (\log_2 r) + \log_2 f(x).
\end{equation}
As a trivial consequence of hypothesis (iii) (which involved the 
natural logarithm) in Theorem \ref{th1.3}, the function 
$x \mapsto \log_2 f(x)$ is convex, and hence the function
$g$ in (\ref{eq3A1}) is convex.
\medskip

     {\bf Claim B.}\ \ {\it The following three statements hold: \zhfb
(1)\ \ $(\log_2 r)x - g(x) \to \infty$ as $x \to \infty$. \zhfb
(2)\ \ For every $t < \log_2 r$, one has that
$g(x) - tx \to \infty$ as $x \to \infty$.} 
\smallskip

   {\it Proof.}\ \ To verify (1), note that by (\ref{eq3A1}) and hypothesis (i) in Theorem \ref{th1.3},
$g(x) - (\log_2 r)x = \log_2 f(x) \to -\infty$ as $x \to \infty$.
Thus (1) holds.

   To verity (2), suppose $t$ is any number such that 
$t < \log_2 r$.
Then $2^t < r$ and hence $2^t/r < 1$.
Hence by hypothesis (ii) in Theorem \ref{th1.3},
$(2^t/r)^x /f(x) \to 0$ as $x \to \infty$.
That is, $2^{tx} / [r^x f(x)] \to 0$ as $x \to \infty$.   
Hence by the first equality in (\ref{eq3A1}),  
$tx - g(x) \to -\infty$ as $x \to \infty$.
Thus (2) holds.
That completes the proof.
\medskip 

     {\bf Step C.}\ \ For convenient later reference, we shall just list here again, with 
substantial redundancy, the properties of the function 
$g: [0,\infty) \to (-\infty, -1]$ in (\ref{eq3A1}) 
[and the related number $\log^2 r$] 
that will be used in the proof.
In this list, item (0) holds trivially by the hypothesis that 
$0 < r \leq 1$;
items (1), (2), and (4) hold by (0) and (\ref{eq3A1}) and the 
hypothesis of Theorem 1.3; 
item (3) repeats the sentence after (\ref{eq3A1});
items (5) and (6) simply repeat the statements in Claim B;
and item (7) holds by (\ref{eq3A1}).   
Here is the list: \smallskip

\noindent (0)\ \ $-\infty < \log_2 r \leq 0$. \zhfb
\noindent (1)\ \ $g$ is negative; in fact $g(x) \leq -1$ for all $x \in [0,\infty)$. \zhfb
\noindent (2)\ \ $g$ is continuous on $[0,\infty)$.  \zhfb
\noindent (3)\ \ $g$ is convex on on $[0,\infty)$.  \zhfb
\noindent (4)\ \ $g$ is strictly decreasing on $[0,\infty)$. \zhfb
\noindent (5)\ \ 
$\lim_{x \to \infty} [(\log_2 r)x - g(x)] = \infty$. \zhfb
\noindent (6)\ \ For every $t < (\log_2 r)$, 
$\lim_{x \to \infty} [g(x) - tx] = \infty$. \zhfb
\noindent (7)\ \ $2^{g(x)} = r^x f(x)$ for all $x \in [0, \infty)$. 
\medskip

   {\bf Step D.}\ \ Suppose $v$ and $y$ are any two real numbers such that
$0 \leq v < y$.

     Define the real number $\zeta_{v,y}$ by 
 \begin{equation}
\label{eq3D1}
     \zeta_{v,y}\ :=\ {{g(y) - g(v)} \over {y-v}}.
\end{equation}   
By Step C(4), $\zeta_{v,y} < 0$.

     Define the (affine) function $C_{v,y}: [v,y] \to \R$ as follows:
\begin{equation}
\label{eq3D2}
    {\rm For\ each}\ x \in [v,y],\ \ \ C_{v,y}(x)\ :=\ g(v) + \zeta_{v,y}(x - v).
\end{equation}
Then $C_{v,y}(v) = g(v)$ and  $C_{v,y}(y) = g(y)$.
This function $C_{v,y}$ simply specifies the 
``chord'' in $\R^2$ from the point
$(v,g(v))$ to the point $(y,v(y))$.
By (\ref{eq3D2}), the inequality $\zeta_{v,y} < 0$ above
(after (\ref{eq3D1})), and then Step C(3)(4), one 
has that
\begin{equation}
\label{eq3D3}
    g(v) \geq C_{v,y}(x) \geq g(x) \geq g(y)\ \ \ 
{\rm for\ all}\ x \in [v,y].
\end{equation}    
Define the nonnegative real number $M_{v,y}$ by           
\begin{equation}
\label{eq3D4}
    M_{v,y} := \sup_{x \in [v,y]}[C_{v,y}(x) - g(x)].
\end{equation}

      {\bf Claim E.}\ \ {\it For any element $v \in [0, \infty)$, there exists $w \in (v, \infty)$ such that $M_{v,w} = 1$.}
\smallskip
  
     {\it Proof.}\ \ Suppose $v \in [0, \infty)$.
     
     By Step C(5),
$[g(v) + (\log_2 r)(x-v) - g(x)]
= [(\log_2 r)x - g(x)] - [(\log_2 r)v - g(v)]
\to \infty$ as $x \to \infty$.
Accordingly, let $z \in (v, \infty)$ be such that  
$[g(v) + (\log_2 r)(z-v) - g(z)] \geq 3$ (say).
Then let $t < \log_2 r$ be such that
$[g(v) + t(z-v) - g(z)] \geq 2$.

By Step C(6), 
$[g(v) + t(x-v) - g(x)] = [tx -g(x)] - [tv - g(v)]
\to -\infty$ as $x \to \infty$.
Applying that and Step C(2) and the last sentence of the preceding paragraph, let $y \in (z, \infty)$ be such that
$g(v) + t(y-v) - g(y) = 0$.
Then $g(v) +t(y-v) = g(y)$.

   From above, one has that $v < z < y$.
By (\ref{eq3D1}) and the last sentence of the preceding 
paragraph just above, $\zeta_{v,y} = t$; and      
hence for each $x \in [v,y]$, 
$C_{v,y}(x) = g(v) + t(x - v)$ by (\ref{eq3D2}).
By the preceding two sentences, followed by the last sentence
of the next-to-last paragraph above,
$C_{v,y}(z) = g(v) + t(z-v) \geq 2 + g(z)$.
Hence by (\ref{eq3D4}),  $M_{v,y} \geq 2$.

   Trivially by (\ref{eq3D4}), (\ref{eq3D3}), and Step C(2),  
$\lim_{u \to v+} M_{v,u} = 0$.
Also, as a consequence of Step C(2)(3), the mapping 
$u \mapsto M_{v,u}$, for $u \in (v,\infty)$, is
(nondecreasing and) continuous.
By the preceding two sentences here and the last sentence of the preceding paragraph, there exists $w \in (v,y)$ such that
$M_{v,w} = 1$.
Thus Claim E holds. 
\medskip

     {\bf Recursion F.}\ \ We shall recursively define an infinite sequence $(y_0, y_1, y_2, \dots)$ of elements of $[0,\infty)$, as follows:

     To start off, define $y_0 := 0$.
   
     Now suppose $n$ is a nonnegative integer, and $y_n$ has already been defined in
$[0,\infty)$.
Applying Claim E, and writing $y_j$ also as $y(j)$ for
typographical convenience, let $y_{n+1}$ be such that 
\begin{equation}
\label{eq3F1}
    y_{n+1} > y_n\ \ \ {\rm and}\ \ \ M_{y(n), y(n+1)} = 1.
\end{equation} 

       That completes Recursion F.
\medskip

     {\bf Claim G.}\ \ {\it $y_n \to \infty$ as $n \to \infty$.}  
\medskip

     {\it Proof.}\ \ Suppose Claim G is false.
Then by the ``first half'' of (\ref{eq3F1}), there is a positive number $z$ such that $y_n \uparrow z$ as $n \to \infty$.
By Step C(2), $g(y_n) \to g(z)$ as $n \to \infty$.
Hence $g(y_{n-1}) - g(y_n) \to 0$ as $n \to \infty$.
Hence by (\ref{eq3D3}) and (\ref{eq3D4}), $M_{y(n-1), y(n)} \to 0$ as $n \to \infty$.
But that contradicts (\ref{eq3F1}).
Thus Claim G holds after all. 
\medskip

     {\bf Step H.}\ \ In this ``step'', we shall simplify some notations, ``extend certain key chords to full lines'', and define some useful reference points $w_n \in (y_{n-1}, y_n)$ for
$n \in \N$.

    Refer to the (\ref{eq3F1}) (its ``first half''), to (\ref{eq3D1}), and to (\ref{eq3D2}).
For each positive integer $n$, first define the real number 
$s_n := \zeta_{y(n-1), y(n)}$ [and keep in mind that 
$s_n < 0$ by the comment after (\ref{eq3D1})], then define the (affine) function $L_n: \R \to \R$ by
\begin{equation}
\label{eq3H1}
     L_n(x) := g(y_{n-1}) + s_n \cdot (x - y_{n-1})\ \ \ {\rm for}\ \ x \in \R,
\end{equation}           
and then define the real number $a_n := L_n(0)$.
From all three of those definitions and some basic algebra, along with (\ref{eq3D2}) and the sentence right after it, one has that for each positive integer $n$,
\begin{align}
\label{eq3H2}
L_n(x) &= a_n + s_nx\ \ \ {\rm for\ all}\ \ x \in \R, \\ 
\label{eq3H3}
L_n(x) &= C_{y(n-1), y(n)}(x)\ \ \ {\rm for\ all}\ \ x \in [y_{n-1}, y_n],
\end{align}
and (hence) in particular, $L_n(y_{n-1}) = g(y_{n-1})$ and $L_n(y_n) = g(y_n)$.

    From (\ref{eq3H3}), (\ref{eq3D4}), and (\ref{eq3F1}),
one has that for each positive integer $n$,
\begin{equation}
\label{eq3H4}
\sup_{x \in [y(n-1), y(n)]} [L_n(x) - g(x)] = M_{y(n-1), y(n)} 
= 1. 
\end{equation}       
Of course for any given $n \in \N$, by Step C(2) and (say) (\ref{eq3H2}), the expression in the main brackets in the left side of (\ref{eq3H4}) is a continuous function of $x$ on
(at least) the closed interval $[y_{n-1}, y_n]$.
Accordingly, for each $n \in \N$, let $w_n \in [y_{n-1}, y_n]$ be such that $L_n(w_n) - g(w_n) = 1$.
For any given $n \in \N$, by the two equations right after (\ref{eq3H3}), one in fact has that
$y_{n-1} < w_n < y_n$.
\medskip

   {\bf Step I.}\ \ Here we shall display for convenient reference seven technical facts that 
will be employed in (the rest of) the proof of Theorem \ref{th1.3}.

     Here are the first two: 
\begin{align}
\label{eq3I1}
&0 = y_0 < w_1 < y_1 < w_2 <  y_2 < w_3 <  y_3 < \dots;\ \ \  {\rm and}\\
\label{eq3I2}
&y_n \to \infty\ \ \ {\rm as}\ \ n \to \infty.
\end{align}
Here (\ref{eq3I1}) comes from Recursion F and (for every positive integer $n$) the final pair of inequalities in Step H.   
Eq.\ (\ref{eq3I2}) is simply Clain G.

   Next, for each positive integer $n$, 
\begin{align}
\label{eq3I3}
   L_n(y_{n-1}) &= g(y_{n-1});\\
\label{eq3I4}
   L_n(w_n) &= g(w_n) + 1;\ \ \  {\rm and}\\
\label{eq3I5}
   L_n(y_n) &= g(y_n).
\end{align}
Eqs.\ (\ref{eq3I3}) and (\ref{eq3I5}) were pointed out right after (\ref{eq3H3}), and eq.\ (\ref{eq3I4}) comes from the 
next to last sentence in Step H.
Also, for each positive integer $n$,          
\begin{equation}
\label{eq3I6}
g(x) \leq L_n(x) \leq g(x) + 1\ \ \ {\rm for\ all}\ \ x \in [y_{n-1}, y_n].
\end{equation}
Here the first inequality comes from (\ref{eq3I1}), 
(\ref{eq3I3}), (\ref{eq3I5}), and Step C(3) (convexity of $g$); and the second inequality come from (\ref{eq3H4}).

     Now in fact (\ref{eq3I1}), (\ref{eq3I3}), (\ref{eq3I4}), 
(\ref{eq3I5}), and Step C(3) (again, convexity of $g$) yield that for each positive integer $n$, one has that
$g(x) < L_n(x)$ for all $x \in (y_{n-1}, y_n)$, and that, again for each $n \in \N$,
\begin{equation}
\label{eq3I7}
    L_n(x) < g(x)\ \ \ {\rm for\ all}\ \ x \in [0, \infty) - [y_{n-1}, y_n].
\end{equation}     

     {\bf Remark J.}\ \ For any given positive integer $n$, the following comments
(1)-(7) hold:

   (1)\ \ By (\ref{eq3I1}), 
$0 \leq y_{n-1} < w_n < y_n < w_{n+1} < y_{n+1}$;
and hence $w_n \in [0, \infty) - [y_n, y_{n+1}]$
and $w_{n+1} \in [0, \infty) - [y_{n-1}, y_n]$.  

   (2)\ \ By (\ref{eq3I3}) and (\ref{eq3I5}), 
$L_{n+1}(y_n) - L_n(y_n) = g(y_n) - g(y_n) = 0$.

   (3)\ \ By (1) above and (\ref{eq3I7}), and then ({\ref{eq3I4}), 
$L_{n+1}(w_n) - L_n(w_n) < g(w_n) - L_n(w_n) = -1$.    

   (4)\ \ By (1), (2), and (3) above, the function
$x \mapsto L_{n+1}(x) - L_n(x)$ for $x \in \R$ [with constant derivative $s_{n+1} - s_n$ --- recall (\ref{eq3H2})] is strictly increasing (and in particular, $s_{n+1} - s_n > 0$).

   (5)\ \ By (1), (3), and (4) above, 
$L_{n+1}(x) - L_n(x) < -1$ for all $x \leq y_{n-1}$ 
(in fact for all $x \leq w_n$).

   (6)\ \ By (1) above and (\ref{eq3I7}), and then ({\ref{eq3I4}),
$L_{n+1}(w_{n+1}) - L_n(w_{n+1}) >  
L_{n+1}(w_{n+1}) - g(w_{n+1}) = 1$.

   (7)\ \ By (1), (6), and (4) above,
$L_{n+1}(x) - L_n(x) > 1$ for all $x \geq y_{n+1}$
(in fact for all $x \geq w_{n+1}$).
\medskip

   {\bf Remark K.}\ \ For any two given positive integers $h$ and $j$ such that $h < j$,
the following comments (1)-(5) hold:

    (1)\ \ For each $x \in \R$, one has the telescoping sum
$ L_j(x) - L_h(x) = \sum_{i=h}^{i = j-1} [L_{i+1}(x) - L_i(x)]$.
 
    (2)\ \ For any $x \geq y_j$ and any 
$i \in \{h, h+1, \dots, j-1\}$, one has that
$x \geq y_{i+1}$ by (\ref{eq3I1}) and hence
$L_{i+1}(x) - L_i(x) > 1$ by Remark J(7).

    (3)\ \ By (1) and (2) above, $L_j(x) - L_h(x) > j-h$ for all 
$x \geq y_j$.

    (4)\ \ For any $x \leq y_{h-1}$ and any 
$i \in \{h, h+1, \dots, j-1\}$, one has that
$x \leq y_{i-1}$ by (\ref{eq3I1}) and hence
$L_{i+1}(x) - L_i(x) < -1$ by Remark J(5).
     
    (5)\ \ By (1) and (4) above, $L_j(x) - L_h(x) < - (j-h)$ 
for all $x \leq y_{h-1}$.
\medskip
  
   {\bf Remark L.}\ \ By (\ref{eq3H2}), (\ref{eq3I1}), (\ref{eq3I3}). and Step C(1), 
\begin{equation}
\label{eq3L1}
   a_1 = L_1(0) = L_1(y_0) = g(y_0) = g(0) \leq -1.
\end{equation}
For any integer $j \geq 2$, by (\ref{eq3H2}) and Remark K(5)
(with $h=1$ and $x = 0$), 
$a_j - a_1 = L_j(0) - L_1(0) < -(j-1)$,
and hence by (\ref{eq3L1}), $a_j < -j$.
Combining that with (\ref{eq3L1}) itself, one has that
\begin{equation}
\label{eq3L2}
   a_n \leq -n\ \ \ {\rm for\ every}\ \ n \in \N.
\end{equation}

 {\bf Remark M.}\ \ For any given positive integer $n$, the
following comments (1)-(6) hold:

     (1)\ \ By (\ref{eq3I6}) and (\ref{eq3I7}), 
$L_n(x) \leq g(x) + 1$ for all $x \in [0, \infty)$.

     (2)\ \ By (\ref{eq3I7}) (for $x < y_n$) and (\ref{eq3I3})
(for $x = y_n$), 
$L_{n+1}(x) \leq g(x)$ for all $x \leq y_n$.

     (3)\ \ By (2) above and Remark K(5),
for any integer $j \geq n+2$ and any $x \leq y_n$,
one has that  
$L_j(x) = L_{n+1}(x) + [L_j(x) - L_{n+1}(x)]
 \leq g(x) - [j - (n+1)]$. 
 
    (4)\ \ If $n \geq 2$, then by (\ref{eq3I7}) 
(for $x > y_{n-1}$) and (\ref{eq3I5}) (for $x = y_{n-1}$),     
$L_{n-1}(x) \leq g(x)$ for all $x \geq y_{n-1}$.

     (5)\ \ If $n \geq 3$, then by (4) above and Remark K(3), for any positive integer $h \leq n-2$ and any $x \geq y_{n-1}$,
$L_h(x) = L_{n-1}(x) - [L_{n-1}(x) - L_h(x)]
\leq g(x) - [(n-1) - h]$.

     (6)\ \ For the given $n \in \N$, comments (1)-(5) apply 
to all $x \in [y_{n-1}, y_n]$, and [see also ({\ref{eq3H2})]
can be repeated together in one convenient display, as follows:
\begin{equation}
\label{eq3M1}
{\rm For\ all}\ x \in [y_{n-1}, y_n]\ 
{\rm and\ all}\ j \in \N,\ \ \
a_j + s_jx\ =\ L_j(x)\ \leq\ g(x) + 1 - |n-j|.
\end{equation} 

     Claim N below will involve sums whose summands are respectively 2 to the power $L_j(x)$ for $j \in \N$.
The left side of (\ref{eq3M1}) was included because it may 
make slightly less cumbersome the notations in both the display in Claim N itself and the subsequent application of Claim N later on.
\medskip

   {\bf Claim N.}\ \ {\it Writing the coefficients $a_n$ and $s_n$ in (\ref{eq3H2}) also
as $a(n)$ and $s(n)$ respectively, one has the following:
For every $x \in [0, \infty)$,}
\begin{equation}
\label{eq3N1}
   \sum_{j=1}^\infty \bigl[2^{a(j)} \cdot (2^{s(j)})^x\bigl]\ \leq\ 6 \cdot r^x f(x).
\end{equation} 

   {\it Proof.}\ \ Suppose $x \in [0, \infty)$.
Applying (\ref{eq3I1}) and (\ref{eq3I2}), let $n \in \N$ be such that $x \in [y_{n-1}, y_n]$.
By Remark M(6), for every positive integer $j$, 
$a_j + s_jx \leq g(x) + 1 - |n-j|$.
Hence by Step C(7), 
\begin{align*}
    \sum_{j=1}^\infty \bigl[2^{a(j)} \cdot (2^{s(j)})^x\bigl]\ &\leq\ \ 
     \sum_{j=1}^\infty 2^{g(x) + 1 - |n-j|}\ =\ 2^{g(x)+1} \sum_{j=1}^\infty 2^{-|n-j|}\\
&\leq\ 2 \cdot r^x f(x) \cdot \sum_{j=-\infty}^\infty 2^{-|n-j|}\ =\ 2 \cdot r^x f(x) \cdot 3.       
\end{align*} 
Thus (\ref{eq3N1}) holds.
That completes the proof.
\medskip        
   
     {\bf Step O.}\ \ For each positive integer $n$, one has that $s_n < 0$ as was noted right before (\ref{eq3H1}), and $s_{n+1} - s_n > 0$ as was noted at the end of
Remark J(4).  
Thus
\begin{equation}
\label{eq3O1}
   s_1 < s_2 < s_3 < \dots < 0.
\end{equation}
Define the number $\eta \in (-\infty, 0]$ by
\begin{equation}
\label{eq3O2}
   \eta\ :=\ \sup_{n \in \N} s_n\  =\  \lim_{n \to \infty} s_n.
\end{equation}
Our next task is to establish the value of this number $\eta$.

   Recall from Step C(5) that $g(x) - (\log_2 r)x \to -\infty$
as $x \to \infty$.
It follows that for any given $n \in \N$, by (\ref{eq3H2})
and (\ref{eq3I7}),
$a_n + s_nx - (\log_2 r)x = L_n(x) - (\log_2 r)x \to -\infty$
as $x \to \infty$,
hence $(s_n - \log_2 r)x \to -\infty$ as $x \to \infty$, 
hence $s_n - \log_2 r < 0$ must hold,
that is, $s_n < \log_2 r$.
Since that holds for all $n \in \N$, it follows from
(\ref{eq3O2}) that $\eta \leq \log_2 r$.

   Next, for any given $n \in \N$ and any given 
$x \in [y_{n-1}, y_n]$, one has that
[$x \geq 0$ by (\ref{eq3I1}) and]
$\eta x \geq s_n x > a_n + s_n x = L_n(x) \geq g(x)$
by (\ref{eq3O2}), (\ref{eq3L2}), (\ref{eq3H2}), and
(\ref{eq3I6}).
Hence $\eta x > g(x)$ for all $x \in [0, \infty)$
by (\ref{eq3I1}) and (\ref{eq3I2}). 
Hence $\eta \geq \log_2 r$ [for otherwise a contradiction would result from Step C(6)].

   By the final sentence of each of the preceding two paragraphs,
$\eta = \log_2 r$.
Hence by (\ref{eq3O2}),
\begin{equation}
\label{eq3O3}
\log_2 r\ =\ \sup_{n \in \N} s_n\  =\  \lim_{n \to \infty} s_n.
\end{equation}  

   \bf Step P.}\ \ Now the Markov chains in Definition \ref{df2.4} will be brought into play.
   
     (1)\ \ First recall from (\ref{eq3L2}) and (\ref{eq3O1}) that for every positive integer $j$,
one has that $a_j \leq -j$ and $s_j < 0$.   
Thus
\begin{equation}
\label{eq3P1}
   {\rm for\ every}\ j \in \N, \ \ \  
0 < 2^{a(j)} \leq 2^{-j} \leq 1/2\ \ \ {\rm and}\ \ \ 0 < 2^{s(j)} < 1
\end{equation}     
where $a_j$ and $s_j$ are (again) written also as $a(j)$ and $s(j)$ for typographical convenience.
\smallskip

      (2)\ \ Now for each positive integer $j$ referring to (\ref{eq3P1}) and 
Definition \ref{df2.4}, let $W^{(j)} := (W^{(j)}_k,\, k \in \Z)$ 
be a strictly stationary Markov chain
that satisfies Condition $\cS(2^{a(j)}, 2^{s(j)})$.
Let these Markov chains $W^{(1)},\,  W^{(2)},\, W^{(3)},\,\dots$ be constructed in such a way that they are independent of each other.
\smallskip

     (3)\ \ By Definition \ref{df2.4}(I)(II), for any given positive integer $j$ and any given
$k \in \Z$, the random variable $W^{(j)}_k$ takes only the values 0 and 1, with probabilities
\begin{equation}
\label{eq3P2}
P(W^{(j)}_k = 0) = 1 - 2^{a(j)}\ \ \ {\rm and}\ \ \ P(W^{(j)}_k = 1) = 2^{a(j)}.
\end{equation}

     (4)\ \ By (\ref{eq3P1}), $\sum_{j=1}^\infty 2^{a(j)} < \infty.$
Hence by (\ref{eq3P2}) and the Borel-Cantelli Lemma, for any given $k \in \Z$, 
\begin{equation}
\label{eq3P3}
P(W^{(j)}_k = 1\ \rm{for\ infinitely\ many}\ j \in \N)\ =\ 0.
\end{equation}     
Just for technical convenience, deleting a set of probability 0 from the probability space 
$\Omega$ if necessary, we assume without loss of generality that for each integer $k$,
the event in the left side of (\ref{eq3P3}) is the empty set.
      
\smallskip
     (5)\ \ In what follows, here (and in Step R below), 
when the notation $W^{(j)}$ 
appears as a subscript, it will be written as $W(j)$ for typographical convenience.
 
By Lemma \ref{lm2.5}, for each positive integer $j$, the following statements hold:      
\begin{align}
\label{eq3P4}
   &{\rm The\ Markov\ chain}\ W^{(j)}\ {\rm is\ reversible.}\\
\label{eq3P5}
   &{\rm For\ each}\ m \in \N,\ \ \rho_{W(j)}(m) = (2^{s(j)})^m.\\
\label{eq3P6}
   &{\rm For\ each}\ m \in \N,\ \ \alpha_{W(j)}(m) \geq (1/2) \cdot (2^{a(j)}) \cdot (2^{s(j)})^m.\\
\label{eq3P7}
   &{\rm For\ each}\ m \in \N,\ \ \beta_{W(j)}(m) \leq 2 \cdot (2^{a(j)}) \cdot (2^{s(j)})^m.
\end{align}

     {\bf Step Q.}\ \ Let $X := (X_k, k \in \Z)$ be the sequence of random variables defined
as follows:  
\begin{equation}
\label{eq3Q1}
{\rm For\ each}\ k \in \Z,\ \ X_k := \sum_{j=1}^\infty 2^{j-1} W^{(j)}_k.
\end{equation}  
By Step P(3)(4), (recall the sentence after (\ref{eq3P3})), for each $k \in \Z$, the sum in
(\ref{eq3Q1}) has, at any given sample point $\omega \in \Omega$, at most finitely many
non-zero terms (and therefore trivially converges).     

     By Step P(3)(4), and (\ref{eq3Q1}), the random variables $X_k, k \in \Z$ take their 
values in the (countably infinite) set $\{0,1,2,\dots\}$ of all nonnegative integers.
In fact ({\ref{eq3Q1}) involves a standard one-to-one correspondence between that
set $\{0,1,2,\dots\}$ and the set of all sequences $(z_1, z_2, z_3, \dots)$ of
elements of $\{0,1\}$ with at most finitely many 1's. 
As a consequence,
\begin{equation}
\label{eq3Q2}
{\rm for\ each}\ k \in \Z,\ \ \sigma(X_k) = \sigma(W^{(j)}_k, j \in \N).
\end{equation}

      Further, by (\ref{eq3P1}) and (\ref{eq3P2}), for any given $k \in \Z$ and any given 
sequence $(z_1, z_2, z_3, \dots) \in \{0,1\}^\N$ with at most finitely many 1's, the numbers     
$P(W^{(j)}_k = z_j),\ j \in \N$ are each an element of the open interval $(0,1)$, and they      
satisfy $\sum_{j=1}^\infty [1 - P(W^{(j)}_k = z_j)] < \infty$, and hence
\begin{equation*}
P\Bigl(\ \bigcap_{j=1}^\infty \{W^{(j)}_k = z_j\}\Bigl)\ =\ \prod_{j=1}^\infty P(W^{(j)}_k = z_j)\
>\ 0,
\end{equation*}
where the equality comes from the second sentence of Step P(2).
It follows from (\ref{eq3Q1}) that for any given $k \in \Z$, and any given 
$m \in \{0,1,2,\dots\}$, $P(X_k = m) > 0$.

     With elementary (if tedious) arguments, using (\ref{eq3Q1}), (\ref{eq3Q2}), Step P(2), and Lemma \ref{lm2.5}(I), one can verify that the random sequence $X$ is strictly stationary, that $X$ is a Markov chain (with state space $\{0,1,2,\dots\}$), and that
$X$ is reversible.
\medskip

     {\bf Step R.}\ \ Recall the last paragraph of Step P and the last sentence of Step Q.
To complete the proof of Theorem \ref{th1.3}, 
what remains is to verify (\ref{eq1.33}).

     Let $m$ be an arbitrary fixed positive integer.   
For convenient reference, taking (\ref{eq1.19}) into account, what remains to prove are the following three items: \zhfb
(1)\ \ $\rho_X(m) = r^m$; \zhfb
(2)\ \ $\alpha_X(m) \geq (1/2) \cdot r^m f(m)$; and \zhfb
(3)\ \ $\beta_X(m) \leq 12 \cdot r^m f(m)$. \zhfb
In the proofs of (1)-(3) given below, when the notation 
$W^{(j)}$ is itself a subscript,
it will be written as $W(j)$ for typographical convenience.
\smallskip

     {\it Proof of (1).}\ \ By eq.\ (\ref{eq3Q2}), Step P(2), Lemma \ref{lm2.6},  Lemma \ref{lm2.5}(III), 
and eq.\ (\ref{eq3O3}), one has that
\begin{equation*}
\rho_X(m)\ =\ \sup_{j \in \N} \rho_{W(j)} (m)\ =\ \sup_{j \in \N} (2^{s(j)})^m\ =\ r^m.
\end{equation*}
Thus (1) holds.
\smallskip

   {\it Proof of (2).}\ \ Applying (\ref{eq3I1}) and (\ref{eq3I2}), let $j \in \N$ be such that
$m \in [y_{j-1}, y_j]$.
Then by (\ref{eq3H2}) and (\ref{eq3I6}), $a_j + s_j m = L_j(m) \geq g(m)$.  
Taking 2 to the powers on the left and right sides and applying Step C(7), one obtains that 
$2^{a(j)} \cdot (2^{s(j)})^m \geq 2^{g(m)} = r^m f(m)$
(where again $a_j$ and $s_j$ are written as $a(j)$ and $s(j)$).
Hence by (\ref{eq3Q2}), followed by Step P(2) and 
Lemma \ref{lm2.5}(III), one has that
$\alpha_X(m) \geq \alpha_{W(j)}(m) 
\geq (1/2) \cdot r^m f(m)$.
Thus (2) holds.
\smallskip

   {\it Proof of (3).}\ \ By eq.\ (\ref{eq3Q2}), Step P(2), 
Lemma 2.6, Lemma \ref{lm2.5}(III), and finally Claim N, 
one has that
\begin{equation*}
\beta_X(m)\ \leq\ \sum_{j=1}^\infty \beta_{W(j)}(m)\ 
\leq\ \sum_{j=1}^\infty\, [2 \cdot 2^{a(j)} \cdot (2^{s(j)})^m]\ \leq\ 12 \cdot r^m f(m). 
\end{equation*}
Thus (3) holds.
That completes the proof of Theorem \ref{th1.3}. 
\medskip


  




\end{document}